\documentclass[11pt]{extarticle}
\usepackage[utf8]{inputenc}
\usepackage[T1,T2A]{fontenc}
\usepackage{geometry}
\usepackage{amsmath,amsfonts,amssymb,amscd,euscript}
\usepackage{amsthm}
\usepackage[dvipsnames]{xcolor}
\usepackage[english]{babel}
\usepackage{graphicx} 
\usepackage{psfrag} 
\usepackage{float}
\usepackage{subfig}
\usepackage{setspace}
\usepackage{bm}
\usepackage{mathabx}
\usepackage{url}
\usepackage{multirow,tabularx}
\usepackage{csquotes}

\usepackage[ruled,vlined]{algorithm2e}

\graphicspath{{Figures/}}

\numberwithin{equation}{section}

\usepackage[markup=underlined]{changes}

\newtheorem{remark}{Remark}

\title{Statistical parameter identification of mixed-mode patterns from a single experimental snapshot}
\author{Alexey Kazarnikov$^1$, Robert Scheichl$^1$, Irving R. Epstein$^2$,  Heikki Haario$^3$\\ and Anna Marciniak-Czochra$^{1}$}

\date{\footnotesize
${}^1$ Institute for Mathematics and Interdisciplinary Center for Scientific Computing (IWR), Heidelberg University, Im Neuenheimer Feld 205, 69120 Heidelberg, Germany\\[0.5ex]
${}^2$ Department of Chemistry and Volen Center for Complex Systems, Brandeis University, Waltham 02454-9110, Massachusetts, United States\\[0.5ex]
${}^3$ Department of Computational Engineering, LUT University, Yliopistonkatu 34, 53850 Lappeenranta, Finland\\[1ex]
}

\begin{document}

\maketitle

\textbf{Abstract:} Parameter identification in pattern formation models from a single experimental snapshot is challenging, as traditional methods often require knowledge of initial conditions or transient dynamics—data that are frequently unavailable in experimental settings. In this study, we extend the recently developed statistical approach, Correlation Integral Likelihood (CIL) method to enable robust parameter identification from a single snapshot of an experimental pattern. Using the chlorite-iodite-malonic acid (CIMA) reaction — a well-studied system that produces Turing patterns — as a test case, we address key experimental challenges such as measurement noise, model-data discrepancies, and the presence of mixed-mode patterns, where different spatial structures (e.g., coexisting stripes and dots) emerge under the same conditions. Numerical experiments demonstrate that our method accurately estimates model parameters, even with incomplete or noisy data. This approach lays the groundwork for future applications in developmental biology, chemical reaction modelling, and other systems with heterogeneous output.

\section{Introduction}

Mechanisms of self-organised pattern formation in biological and ecological systems have been a central focus of experimental and theoretical research for several decades \cite{VanHarten1995,Murray1993,Sun2021,Tomlin2007}. Various spatial patterns and morphologies have been mathematically modelled \cite{Gierer1972,Kondo1995,Marciniak2006,Torii2012}. Despite being developed under different hypotheses, these models share common mathematical features and are primarily based on a limited range of pattern generation paradigms \cite{Sun2021,Veerman2021}. The most widely used classical Turing approach is based on the concept of spreading morphogens \cite{Turing1952}. Nonlinear interactions facilitated by contrasting diffusivities may yield destabilisation of homogeneous constant states and development of spatially heterogeneous Turing patterns.  Furthermore, the formation of patterns can emerge from the coupling of diffusive and spatially localised non-diffusive components, a phenomenon that is prevalent in both cellular systems and ecosystems  \cite{Halatek2018,Harting2017}. 
Recent studies have also identified a range of mechanisms that are based on mechanical interactions
\cite{Mercker2018,Mercker2013,weevers2024}, nonlocal feedbacks \cite{Martinez-Garcia2014} or are induced by spatial heterogeneity of the underlying environment \cite{Cantrell2004,He2024,Sun2021}. 

To validate the mechanism behind a mathematical model, it should be verified against experimental observations. However, several factors complicate this process in the context of pattern formation. First, in many experimental situations, it is only possible to observe the stationary regime of the pattern formation process, without detailed knowledge of the transient behaviour or the initial state. Due to instabilities in the underlying nonlinear dynamics, small changes in the initial conditions lead to different patterns, which creates challenges for parameter identification. Additionally, data are often available only in a normalised form, providing information about the topology of the pattern but not quantitative concentration levels. It may also occur that some pattern components are unavailable due to the impossibility of measuring all involved variables simultaneously, or when the model mechanism includes hypothetical quantities. Finally, the observed pattern may not be stationary but can still be in transition to a steady state.   

The problem of parameter identification from a single snapshot of a pattern has been the subject of extensive research in the context of Turing models, which can be regarded as a prototype class of pattern formation models. The applicability of a particular strategy is contingent on the assumptions made regarding the information encoded in the pattern. Some approaches, based on classical optimal control theory \cite{Chang2023a, Ding2024,Gao2022,Garvie2007,Garvie2014,  Li2024,Miao2024,Zhu2023} or Bayesian methods \cite{Campillo-Funollet2019}, require knowledge of the initial data that give rise to the specific pattern. When this information is unavailable, but all components of the Turing pattern are known, one can employ the stationarity of the pattern, by substituting the steady-state solution into the model equations and minimising the residual with respect to model parameters. This leads to approaches based on optimal control theory \cite{Garvie2010}, which simulate the model using the stationary pattern as initial data and minimise the drift caused by incorrect parameter values, or least squares methods \cite{Chang2023,Glasner2021,Sgura2019}, which substitute the pattern into the model equation and minimise the norm of the discrepancy between data and model with respect to model parameters. When some components of the stationary solution are unavailable or only the pattern topology is known, the problem becomes more challenging. The work by \cite{Sgura2019} avoids this issue by applying the least squares in two stages: first, conducting experiments with synthetic data and known initial conditions to cluster the parameter space into regions corresponding to different types of patterns, and then fitting the model to the min-max scaled experimental data of one model component separately within each cluster. Papers \cite{Garvie2014,Matas-Gil2024} address the data limitation by assuming some \enquote{realistic} values of the component amplitudes, making it possible to use optimal control theory and physically-informed neural networks (the least squares approach), respectively. Another robust way to handle the data limitation problem involves using convolutional neural networks to learn the mapping from pattern images to model parameters using either raw data \cite{He2024, Zhu2022} or invariant representations that capture the spatial structure of Turing patterns \cite{Schnoerr2023}. However, the amount of required training data grows exponentially with the number of model parameters, limiting the applicability of these approaches to situations where the number of parameters is small \cite{Kazarnikov2023}. 

In this paper, we extend our previously developed statistical approach \cite{Kazarnikov2020,Kazarnikov2023a,Kazarnikov2023} to deal with real data. We focus on the example of pattern formation in the chlorite-iodite-malonic acid (CIMA) chemical reaction, a classical experimental system capable of producing Turing patterns. Working with experimental data introduces challenges that are not present when working with synthetic data only, such as measurement noise, a possible discrepancy between model and data, or the lack of a reasonable initial guess for the model parameters. Furthermore, a key feature of the system under consideration here is its capacity to generate mixed mode patterns that incorporate multiple characteristic length scales. This is evidenced by the coexistence of shapes such as stripes and spots. Such patterns have been observed experimentally in the CIMA reaction \cite{Rudovics1996}, but also shown in simulations in different Turing-type models \cite{Bozzini2015,Ma2019,Woolley2021}. To the best of our knowledge, parameter estimation for such cases has not been studied so far. Our previously developed Synthetic Correlation Integral Likelihood (SCIL) method cannot be directly applied, since quantifying distances between such heterogeneous patterns 
does not lead to a well-defined single Gaussian likelihood.
In this paper, we present a novel approach to deal with this problem by quantifying multiple distributions.   

Random simulation outcomes {are} not unique to Turing models. They appear in chaotic dynamical systems \cite{Haario2015, Springer2019, Springer2021}. Other classes of pattern formation models, such as neuron growth models \cite{Duswald2024}, phase separation equations \cite{Burger2013}, or cellular automata \cite{Kazarnikov2023a} show similar behaviour. The extension of our approach presented here applies equally well to any of those application fields, where the coexistence of qualitatively different behaviour types emerges \cite{Bouali2012}.

The paper is organised as follows. Section \ref{sec:introduction2} gives a general introduction to statistical parameter identification from pattern data with the CIL approach. Section \ref{sec:chemical_experiment} introduces the chemical experiment that was used to produce the pattern data and discusses the main challenges introduced by experimental data. The main concepts of the CIL algorithm and its extension to the case of mixed mode patterns are described in Section~\ref{sec:cil_approach}. In Section~\ref{sec:parameter_identification}, the numerical experiments are presented and the validation of our results through chemical derivation is discussed.

\section{Robust parameter estimation from a single experimental snapshot\label{sec:introduction2}}

Although many papers address the problem of parameter identification from a single pattern, most of them use synthetic data produced by the underlying model to validate their methods. Very few works consider the more challenging case of fitting a reaction-diffusion model to real experimental data, where only one pattern component is typically available. One of the methods that avoids this problem is the two-stage least squares approach \cite{Sgura2019}. This is achieved by first separating the parameter space into regions corresponding to different pattern types, and then optimising the cost function within the pre-selected cluster. While this approach is successful for two free parameters, it involves manual clustering of the parameter space, which can be problematic when dealing with many parameters. It also requires {\it a priori} knowledge of the parameter range where the correct parameter values are located. The study by \cite{Garvie2014} tackles the problem of partial information by assuming some \enquote{realistic} values of the concentration levels, allowing them to fit the Gierer-Meinhardt model to fish skin patterns. However, in \cite{Garvie2014} the model parameters are assumed to depend on space and time, which in some sense \enquote{encodes} the final pattern into the model parameters, rather than finding the correct regime of the underlying model. Similar assumptions about the experimental data from a chemical reaction are made in \cite{Matas-Gil2024}. While this results in reasonable estimates for synthetic patterns, it does not work accurately for real experimental data.


Correlation Integral Likelihood (CIL) is a statistical approach that can be used to identify parameters of pattern formation models using pattern data. Initially developed for chaotic systems \cite{Haario2015, Springer2019}, it was later extended to pattern formation models and successfully applied to classical Turing models \cite{Kazarnikov2020}, to rumour propagation models both in continuous spatial domains and on networks \cite{Li2024,Zhu2022}, and to cellular automata \cite{Kazarnikov2023a}. In all those applications, it was assumed that multiple observations of the pattern (at least 50) are available. By introducing a synthetic version of the CIL (SCIL) the approach was later extended to situations with limited data \cite{Kazarnikov2023}. This extension enabled the identification of parameters from a single snapshot of a pattern, provided the synthetic data are  well-controlled and the pattern is regular. The approach is not applicable to mixed mode patterns. 

Our approach is innovative in its use of spontaneous pattern formation, quantifying the variability of a family of patterns for given parameter values and leveraging this information to compare model predictions with actual data. This statistical formulation allows us to work without the exact knowledge of the initial data that produced a particular pattern, a requirement of most other methods, and instead use a weaker assumption about the distribution of initial data. Furthermore, it is not assumed that the pattern data contains exhaustive information about all model variables. Consequently, we can operate in situations where only some model components can be experimentally observed, even in scaled or normalised form, which is a common experimental situation in developmental biology.

\begin{figure}
    \centering
    \includegraphics[width=0.99\textwidth]{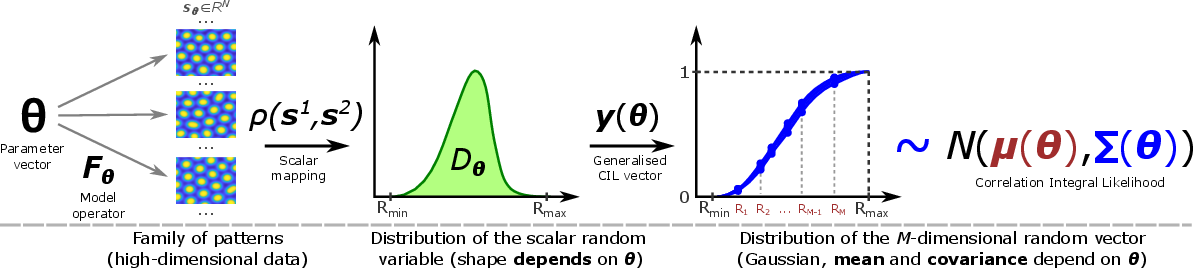}
    \caption{The main idea of the Correlation Integral Likelihood. To characterise a family of high-dimensional patterns, coming from an abstract pattern formation model for fixed values of parameters, a scalar feature (a distance $\rho$ between two arbitrary patterns) is used to map a family of random patterns to scalar values.  A number of eCDF vectors are computed, each for a sufficiently large sample of scalars. The mean and covariance of the vectors are estimated, and the multidimensional Gaussianity is numerically verified. This distribution is used to quantify the statistical distance between a pattern formation model and pattern data.}
    \label{fig:cil_scheme}
\end{figure}

The main steps of the CIL approach are illustrated in Figure~\ref{fig:cil_scheme}. The method is based on a generalisation of the central limit theorem (CLT). According to the CLT, the average of a random variable with finite first and second moments converges to a normal distribution. Provided sufficient repeated measurements are available this allows the use of a Gaussian likelihood for parameter estimation. However, in the non-Gaussian case, when the higher moments do not vanish, the cumulative distribution function (CDF) provides more accurate statistics. For a set of scalar-valued data, the empirical CDF (eCDF) directly approximates the underlying statistical distribution. In probability theory, Donsker's theorem \cite{Donsker1952} represents a significant extension of the CLT. It states that the cumulative distribution function (CDF) of independent and identically distributed (i.i.d.) scalar samples asymptotically approaches a Brownian bridge.
We use the CDFs in an approximate form for finite data.  For i.i.d. scalar data, the eCDF vector
at  $M$ selected bin values becomes a $M$-dimensional Gaussian random vector.
 In many applications, the data is not i.i.d., but under suitable conditions on the dependence of the data, Gaussianity still holds \cite{Borovkova2001,Neumeyer2004}. However, in the case of pattern formation models the initial values are always randomised, and the resulting scalar values are indeed i.i.d. 
It is important to note though that the approach is only approximate for finite sample size.  Due to the limited range of the eCDF vectors in the interval $(0,1)$, Gaussianity only holds for bin values that are not located in the tails.  In the context of numerical applications, standard scalar normality tests can be applied to verify the normality of the bin values in use, and the $M$-dimensional  $\chi^2$ test can be used to confirm the  Gaussianity of the eCDF vectors.

For scalar data, the CIL approach is illustrated by the second arrow in Figure~\ref{fig:cil_scheme}: Given a set of samples from the (potentially complicated, scalar) distribution, we compute the respective eCDF vector at a finite number of bin values. This process is repeated sufficiently many times to obtain good estimate of the mean and covariance matrix of the bin values, which then uniquely define the multivariate Gaussian likelihood of the eCDF vector.  

However, our data are inherently high dimensional. Thus, first a scalar-valued mapping must be employed that maps the high-dimensional distribution to a scalar one, which is then used subsequently to construct the eCDF vectors; see \cite{Haario2015,Kazarnikov2020,Springer2019} for earlier examples. This step is illustrated by the first arrow in Figure~\ref{fig:cil_scheme}. Mapping a family of high-dimensional patterns to a scalar random variable inevitably leads to a loss of information. This can be mitigated by mapping the data separately to several scalar-valued distributions and computing a number of eCDF vectors that are concatenated to form a higher dimensional \enquote{feature vector}. Inspired by a recent analytical study~\cite{Kowall2022}, we use different Lebesgue and Sobolev norms that characterise differences between the pattern values, as well as their spatial gradients. This allows us to better quantify the variability of patterns and increase the accuracy of the parameter identification. It is important to note that a concatenation of Gaussian vectors remains Gaussian.  All technical details about the derivation of the CIL approach and its numerical implementation are given in Section~\ref{sec:cil_approach}.

Our work employs distances between pattern data to characterise pattern families. This feature is highly versatile and can be defined for various types of biological and chemical data, such as concentration patterns, cellular automata patterns, and DNA sequences. However, the approach is not limited to the use of distances. When some problem-specific, meaningful scalar-valued characteristics are available, they can also be employed. This was demonstrated, e.g.,  in \cite{Kazarnikov2023a} by using particle sizes for cellular automata, or residuals for synchronization in \cite{Shah2023}. While we here focus on steady-state cases, the approach is also applicable to transient data or to patterns coming from stochastic models, such as agent-based models or cellular automata \cite{Kazarnikov2023a}.



\section{From simulated to experimental pattern data \label{sec:chemical_experiment}}


\subsection{Choice of the prototype model and the corresponding experimental setting}

In this section, we present the specific pattern formation model used later as a test case for parameter identification from real experimental data.
We focus on a reaction-diffusion model describing chlorite-iodide-malonic acid (CIMA) reactions. This particular experimental scenario has been chosen because it closely reflects the challenges commonly encountered in modelling pattern formation. The difficulties encountered relate to the observation of solely stationary patterns without initial and transition processes, the lack of information on some model components, observations limited to the shape of the spatial structure without actual concentration values, and discrepancies between the idealised model mechanism and the actual underlying biological processes. The employment of a system in which the chemical mechanism of the reaction is thoroughly understood may also enable the verification of the numerical results against the original experimental process.
 
 The CIMA reaction is an example of Turing pattern formation in a controlled experimental setting. The Turing theory is based on a bifurcation from a stationary spatially homogeneous distribution driven by contrasting diffusivity in a reaction-diffusion system of at least two components \cite{Murray1993,Turing1952}. Despite its extensive use in various biological applications, the Turing framework is challenging to identify at the molecular level. In this context, the CIMA reaction can serve as a valuable experimental proof of concept, as it provided the first experimental identification of a Turing pattern \cite{Rudovics1996}. In this reaction, the interactions of iodide and chlorite follow the so-called activation-inhibition loop, which is a prerequisite for the Turing instability \cite{Veerman2021}. The required difference in diffusivity is achieved by reducing the mobility of iodide through complexation with an additional immobile (or nearly immobile) species. Furthermore, the resulting complex undergoes a colour change in accordance with the concentration of the activator, thereby rendering the pattern discernible. The scheme of this chemical experiment is shown in Figure~\ref{fig:chemical_experiment}.
Open spatial reactors are designed to overcome the limitation of observing only transient processes in closed reactors. They maintain the reaction systems at a controlled distance from equilibrium. A two-sided open spatial reactor consists of a gel block in contact with two separated reservoirs (A and B). The small pores in the gel suppress fluid motion, ensuring that pattern formation occurs solely due to chemical reactions and diffusion. The gel block is confined between two wide, thin, uniform, porous, and transparent glass plates. These plates allow the diffusion of chemicals and visual observation of the pattern formed between them. Reservoirs A and B contain mixtures of chemicals that are inert individually but react when combined. By using pumps, these reservoirs function as Continuously Fed Stirred Tank Reactors (CSTR), where contents are mixed thoroughly to achieve spatially uniform concentrations. These mixtures are then pumped past the outer sides of the porous plates, allowing chemicals to diffuse through the plates into the gel. Reaction products subsequently diffuse out and are swept away.

\subsection{Mathematical model}
\label{sec:mathmodel}


The Lengyel-Epstein model of the CIMA reaction is based on a simplified mechanism consisting of the following four reactions:
\begin{equation}
    \begin{array}{c}
        \text{MA} + \text{I}_2 \rightarrow \text{IMA} + \text{I}^- + \text{H}^+;
        \quad
        \text{ClO}_2 + \text{I}^- \rightarrow \text{ClO}_2^- + 1/2\,\text{I}_2; \\
        \text{ClO}_2^- + 4\text{I}^- + 4\text{H}^+ \rightarrow 2\text{I}_2 + \text{Cl}^- + 2\text{H}_2\text{O};
        \quad
        \text{S} + \text{I}_2 + \text{I}^- \rightleftharpoons \text{SI}_3^-,
    \end{array}
    \label{eq:chemical_equations}
\end{equation}
where $\text{MA}$, $\text{I}_2$, $\text{ClO}_2$, $\text{ClO}_2^-$, and $\text{I}^-$ are independent variables, $[\text{H}^+]$ is assumed to be constant, $\text{S}$ is a complexing agent, and $\text{Cl}^-$ and $\text{IMA}$ are inert products. The reaction rates for equations \eqref{eq:chemical_equations} are defined by experimentally estimated nonlinear dependencies on the concentrations \cite{ Kadar1995,Lengyel1990}:
\begin{equation}
\begin{array}{c}
r_1 = \dfrac{k_{1a}[\text{MA}][\text{I}_2]}{k_{1b} + [\text{I}_2]};
\,
r_2 = k_2[\text{ClO}_2][\text{I}^-]; \\
r_3 = k_{3a}[\text{ClO}_2^-][\text{I}^-][\text{H}^+] + \dfrac{k_{3b}[\text{ClO}_2^-][\text{I}_2][\text{I}^-]}{\alpha + [\text{I}^-]^2};
\,
r_4 = k_{4a}[\text{S}][\text{I}^-][\text{I}_2] - k_{4b}[\text{SI}_3^-],
\end{array}
\label{eq:reaction_rates}
\end{equation}
where terms in brackets denote mixed concentration values, scalars $k_{1a}$, $k_{1b}$, $k_2$, $k_{3a}$, $k_{3b}$, $k_{4a}$, and $k_{4b}$ are kinetic parameters. The parameter $\alpha$, present in the expression for $r_3$, is an artificial (\textit{ad hoc}) constant. Experimental studies have demonstrated that the second term in the rate $r_3$ in \eqref{eq:reaction_rates} can be accurately approximated by the functional relationship $k_{3b}[\text{ClO}_2^-][\text{I}_2]/[\text{I}^-]$ over a wide range of concentration values \cite{Kern1965}. However, as $\text{I}^- \rightarrow 0$, $r_3 \rightarrow +\infty$. To address this, it was proposed \cite{Weitz1984} to replace the original term $1/[\text{I}^-]$ with $[\text{I}^-]/(\alpha + [\text{I}^-]^2)$. This modification results in a finite rate $r_3$ at $[\text{I}^-] = 0$, which is not chemically realistic. Consequently, $\sqrt{\alpha}$ can be considered a cut-off concentration, below which the rate expression is invalid. In this study, we fix $\alpha = 10^{-15}$, following the approach in \cite{Lengyel1990,Rudovics1999}.

By applying the reaction rates \eqref{eq:reaction_rates} to reaction equations \eqref{eq:chemical_equations} and taking into account the diffusivity of the mobile components, we arrive at the following system of reaction-diffusion equations
\begin{equation}
\begin{array}{rlrl}\vspace{2mm}
\dfrac{\partial [\text{ClO}_2]}{\partial t} &= -r_2 + D_{\text{ClO}_2} \Delta [\text{ClO}_2],
&
\dfrac{\partial [\text{ClO}_2^-]}{\partial t} &= r_2 - r_3 + D_{\text{ClO}_2^-} \Delta [\text{ClO}_2^-], \\\vspace{2mm}
\dfrac{\partial [\text{I}_2]}{\partial t} &= -r_1 + \dfrac{1}{2}r_2 + 2r_3 - r_4 + D_{\text{I}_2} \Delta [\text{I}_2],
&
\dfrac{\partial [\text{MA}]}{\partial t} &= -r_1 + D_{\text{MA}} \Delta [\text{MA}],
\\
\dfrac{\partial [\text{I}^-]}{\partial t} &= r_1 - r_2 - 4r_3 - r_4 + D_{\text{I}^-} \Delta [\text{I}^-],
&
\dfrac{\partial [\text{SI}_3^-]}{\partial t} &= r_4,
\quad
\dfrac{\partial [\text{S}]}{\partial t} = -r_4,
\end{array}
\label{eq:rd_system}
\end{equation}
posed in a flat bounded domain $\Omega = (0, \ell)^2$ with Neumann (zero-flux) boundary conditions. Here, $\Delta = \frac{\partial^2}{\partial x_1^2} + \frac{\partial^2}{\partial x_2^2}$, $\bm{x} = (x_1, x_2) \in \Omega$, $\ell > 0$ denotes the size of the physical domain, and $t > 0$ is dimensional time. Because changes in the concentrations of $[\text{ClO}_2^-]$ and $[\text{I}^-]$ occur much faster than in other reaction components, the concentrations of all other components can be approximated as constants. Additionally, it is assumed that species S is either bound to a gel matrix or is so large that the diffusivities of S and $\text{SI}_3^-$ are
assumed to be zero \cite{Lengyel1991,Lengyel1992}. With these assumptions, and after applying the variable change
\begin{equation}
\begin{array}{c}
    v = \dfrac{[I^-]}{\sqrt{\alpha}},
    \,
    w = \left( \dfrac{k_{3b}[I_2]}{\alpha k_2 [ClO_2]} \right)[ClO_2^-],
    \,
    d = \dfrac{D_{ClO_2^-}}{D_{I^-}}, \\ \\
    a = \left( \dfrac{k_{1a}[MA]}{\sqrt{\alpha} k_2 [ClO_2]} \right) \left( \dfrac{[I_2]}{k_{1b} + [I_2]} \right),
    \,
    b = \dfrac{k_{3b} [I_2]}{\sqrt{\alpha} k_2 [ClO_2]},
    \\ \\
    \sigma = 1 + \dfrac{k_4}{k_{-4}}[S]_0[I_2],
    \,
    t' = k_2[ClO_2]t,
    \,
    x' = \left( \dfrac{k_2 [ClO_2]}{D_{I^-}} \right)^{\frac{1}{2}}x,
\end{array}
\label{eq:variable_change}
\end{equation}
we arrive at the Lengyel-Epstein model \cite{Lengyel1991, Lengyel1992}
\begin{equation}
   \dfrac{\partial v}{\partial t} = \dfrac{1}{\sigma}\left(a - v - 4\dfrac{vw}{1+v^2} + \Delta v \right);
				\dfrac{\partial w}{\partial t} = b \left( v - \dfrac{vw}{1+v^2} \right) + d \Delta w,
    \label{eq:LEmodel}
\end{equation}
where we dropped the primes from $t$ and $x$ to simplify the notation. After rescaling, the system is posed in a dimensionless spatial domain $\Omega_0 = (0, L)^2$ with $L = \left( \frac{k_2 [ClO_2]}{D_{I^-}} \right)^{\frac{1}{2}}\ell$.


\subsection{Experimental setup and the resulting chemical pattern data}

\begin{figure}
    \centering
    \includegraphics[width=0.99\textwidth]{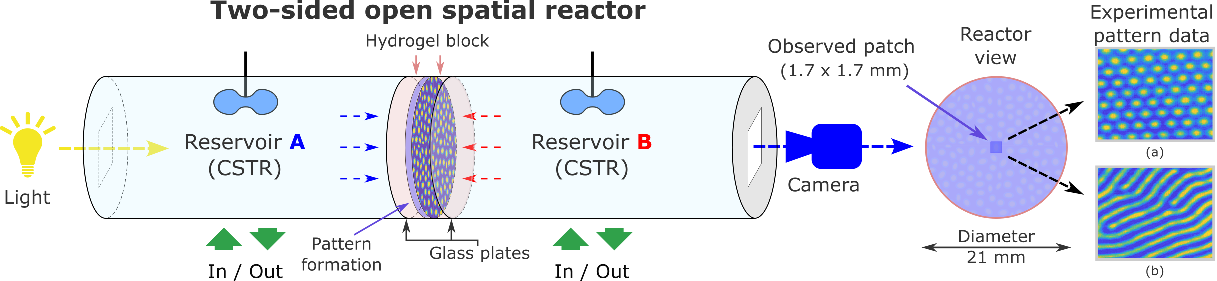}
    \caption{The observation of a chemical Turing pattern in the chlorite-iodite-malonic acid (CIMA) reaction performed in a two-sided open spatial reactor as described in  \cite{Rudovics1996}. The reactor consists of two reservoirs A and B (CSTRs), separated by a hydrogel block, that suppresses convection effects. Each reservoir contains chemical mixtures that are inert individually but react when combined. The hydrogel is loaded with a species of reduced mobility, which slows down the diffusion of the activator and changes colour depending on its concentration, making the pattern visible. A camera captures a digital image of the chemical pattern through a glass window on one side of the reactor.
    }
    \label{fig:chemical_experiment}
\end{figure}

The pattern data used in this study has been sourced from previously published chemical experiments~\cite{Rudovics1996}. Here, we recall the main properties of the experimental conditions, which are necessary to properly compare the model \eqref{eq:LEmodel} with the experimental data. The reaction was performed in a two-sided open spatial reactor using a gel matrix made of agarose, loaded with polyvinyl alcohol (PVA), which served both as a colour indicator and a complexing agent $[\text{S}]$ to reduce mobility. The PVA formed a reddish-purple complex in the presence of iodine and iodide. To maintain consistency, PVA was also fed into the reservoirs in the same concentration to prevent its diffusion from the gel. Chlorite and iodate in basic solutions were fed into Reservoir A, while iodide and malonic acid were fed into the other side in an acetic acid solution. Different patterns were established by varying the concentrations of malonic acid and potassium iodide. Observations were made with a video camera 
with a macro lens. Images were sent to a frame grabber and contrasts were subsequently enhanced. Examples of the patterns used in this paper are shown in Figure~\ref{fig:chemical_experiment} in MATLAB colour scheme.

Calibrating a pattern formation model using the data coming from a real chemical system introduces new challenges not present in the simulated data in our previous studies \cite{Kazarnikov2020,Kazarnikov2023a,Kazarnikov2023}, namely:

\begin{enumerate}
    \item \textbf{Localised observation area.} Unlike in our prior studies, which considered information across the entire computational domain, we now observe exclusively a single square patch in the middle of the reactor (see Figure~\ref{fig:chemical_experiment}). This localised view must be accurately reflected in the numerical simulations of the model. For more details, see Section~\ref{sec:parameter_identification}. 
    \item \textbf{Single component information.} While our earlier works employed both activator and inhibitor components, we now possess information solely about the pattern created by the activator. Consequently, our data are incomplete, adding further complexity to the parameter identification process.
    \item \textbf{Model-data deviation.} Data from real chemical experiments naturally contain imperfections and deviations from the idealised mathematical model. These deviations arise both from incomplete knowledge of the chemical reactions and from the physical characteristics of the measurement devices used to obtain the images. This discrepancy between model and data presents a challenge for parameter identification, but also reflects the experimental situations often encountered in biological systems.
    \item \textbf{Measurement noise and image data.} While simulated data contain actual concentration values of the activator and inhibitor, real data are obtained by photographing the reactor using a digital camera. The resulting greyscale image has pixel values ranging from 0 to 1 (or 0 to 255), which mathematically corresponds to scaling the activator concentration to the interval $(0,1)$, known as a min-max normalisation.
\end{enumerate}

Although such challenges are quite common for any experimental data, they render many of the aforementioned approaches for parameter identification inapplicable. An additional challenge comes from the fact that this particular experimental system  produces mixed mode patterns, exhibiting dots and stripes, for the same parameter value under certain experimental conditions \cite{Rudovics1996}. The CIL approach, as developed in our previous works, cannot be directly applied to this situation. How to modify the CIL approach to overcome this problem is shown in the following section.


\section{Extending the CIL approach to mixed mode patterns \label{sec:cil_approach}}

In this section, we focus on the statistical approach for parameter identification by pattern data and introduce its extension to the new problem of mixed mode patterns. 
 For convenience, we begin with a concise overview of the implementation of CIL for synthetic data, both for large data sets \cite{Kazarnikov2020}, and for the case of a single experimental pattern \cite{Kazarnikov2023}. We then present its extension needed to problems with real data and mixed mode patterns.


\subsection{Formulation of the method}

Let us consider an abstract pattern formation model $\bm{s} = \bm{F}_{\bm{\theta}} (\bm{x})$ on a bounded domain in $\mathbb{R}^2$,
where $\bm{\theta} \in \mathbb{R}^p$ is a vector of control parameters, $\bm{x} \in X$
represents the (spatially distributed) initial data, $\bm{s} \in Y$
is the output pattern, and $\bm{F}_{\bm{\theta}} : X  \rightarrow Y$ for fixed $\bm{\theta} \in \mathbb{R}^p$
is the model operator that maps the initial condition to the final pattern, such as the Lengyel-Epstein model in Section \ref{sec:mathmodel}. In practice, the model operator is assumed to be discretised, e.g., on some uniform grid or into \enquote{pixels}, such that $Y = \mathbb{R}^N$ for some $N \in \mathbb{N}$.
In the models of pattern formation, we observe that patterns $\bm{s}$ change not only due to variations in the model parameters $\bm{\theta}$, but also for fixed model parameters, e.g., due to small variations in the initial condition. 
This is typical for Turing models, where initial conditions are usually defined by perturbing a homogeneous steady state with small random noise, leading to different stationary patterns observed in two-dimensional domains. Thus, we assume that the initial data $\bm{x}$ is a realisation of a random variable $\bm{X}(\omega)$, uniformly distributed around some homogeneous steady state. This induces a probability distribution on the patterns $\bm{s}$, and we denote by $\bm{S}_{\bm{\theta}}(\omega)$ the family of patterns obtained by repeatedly randomising initial data $\bm{x}$ for a fixed parameter $\bm{\theta}$. Provided $\bm{F}_{\bm{\theta}}$ is measurable, it induces a distribution for the high-dimensional, multivariate random variable $\bm{S}_{\bm{\theta}}(\omega) \in \mathbb{R}^N$, which is inherited from the distribution of the initial data $\bm{X}(\omega)$.

Since $\bm{S}_{\bm{\theta}}(\omega)$ is typically very high-dimensional, the core idea of the CIL approach is to map it to a scalar random variable $D_{\bm{\theta}}(\omega)$.  Let  $\rho : \mathbb{R}^N \times \mathbb{R}^N \rightarrow \mathbb{R}_{\geq 0}$ represent a distance between two patterns, based on a suitable norm or metric. This mapping is utilised to compare different pattern families during parameter identification, as discussed in Section~\ref{sec:introduction2}. To characterise the family of patterns generated by the model $\bm{F}_{\bm{\theta}}$ for a fixed parameter vector $\bm{\theta}$, we introduce the scalar random variable $D_{\bm{\theta}} = \rho(\bm{S}^{1}_{\bm{\theta}}, \bm{S}^{2}_{\bm{\theta}})$, where $\bm{S}^{i}_{\bm{\theta}}, i=1,2,$ are two copies of the random variable $\bm{S}_{\bm{\theta}}$,
induced by two i.i.d. copies $\bm{X}^{1}, \bm{X}^{2}$ of the initial values. 
Since patterns are typically bounded, we assume that $D_{\bm{\theta}}$ is also bounded and define 
$R^\text{max}_{\bm{\theta}}  = \text{ess}\,\sup_{\bm{s}_{\bm{\theta}}^1, \bm{s}_{\bm{\theta}}^2} \; \rho(\bm{s}^1_{\bm{\theta}}, \bm{s}^2_{\bm{\theta}})$ where $\bm{s}^1_{\bm{\theta}}, \bm{s}^2_{\bm{\theta}}$ are arbitrary realisations of $\bm{S}_{\bm{\theta}}$.

Consider a set of i.i.d. samples $\{ d^1_{\bm{\theta}}, \ldots, d^n_{\bm{\theta}} \}$ of the random variable $D_{\bm{\theta}}$ obtained by using independent copies of the initial conditions in $n$ runs of the model. We define by $\bm{y}(\bm{\theta}) \in \mathbb{R}^M$ the empirical cumulative distribution function (eCDF) of the distances $d^1_{\bm{\theta}}, \ldots, d^n_{\bm{\theta}}$, which approximates the cumulative distribution function $F_{D_{\bm{\theta}}}$ of  $D_{\bm{\theta}}$ with  a finite number of observations. The components $y_k(\bm{\theta})$ are given by the formula:
\begin{equation}
y_k(\bm{\theta}) = \dfrac{1}{n} \sum \limits_{i = 1}^n \# (d^i_{\bm{\theta}} < R_k), \quad k = 1,\ldots,M,
\label{eq:cil_vector}
\end{equation}
where $0 
< R_1 < \ldots < R_M < R^\text{max}_{\bm{\theta}}$ are fixed threshold values for the bins. 
Due to the Central Limit Theorem (CLT), $\bm{y}(\bm{\theta}) \in \mathbb{R}^N$ asymptotically follows a multivariate Gaussian distribution as $n\rightarrow +\infty$. The probability of a random variable $D^i_{\bm{\theta}}$, $i = 1,\ldots,n$ falling below or above the threshold value $R_k$, $k=1,\ldots,M$, is Bernoulli distributed with expected value $p_k = F_{D_{\bm{\theta}}}(R_k)$ and variance $p_k(1-p_k)$. Thus, the averages in \eqref{eq:cil_vector} converge to a Gaussian distribution as $n \rightarrow +\infty$ by the CLT. 
 Due to the close connection with the construction of the correlation dimension \cite{Haario2015}, the vector $\bm{y}(\bm{\theta}) \in \mathbb{R}^M$ is called the \emph{generalised correlation integral vector},
and the respective distribution the \textit{Correlation Integral Likelihood (CIL)}. The mean and the covariance matrix of the Gaussian distribution are denoted by $\bm{\mu}(\bm{\theta})$ and $\bm{\Sigma}(\bm{\theta})$.  In practice, empirical estimates for 
 $\bm{\mu}(\bm{\theta})$ and $\bm{\Sigma}(\bm{\theta})$ are obtained by repeatedly evaluating the vector 
 $\bm{y}(\bm{\theta}) $ and computing the mean and covariance of the emerging matrix. In real cases, we always deal with a limited number of samples, and need to use resampling to obtain the necessary estimates. As a result, the data independence is lost, but the Gaussianity still holds. The details are discussed in the next subsections.

\subsection{Parameter identification by pattern data}

The problem of parameter identification using pattern data can be formulated as follows: given a finite set $S_\text{data}$ of $N_\text{data}$ patterns, find the parameter vector $\bm{\theta}_0$ that minimises the discrepancy between $S_\text{data}$ and the family of patterns described by a set of realisations of the random variable $\bm{S}_{\bm{\theta}}$ obtained from repeated simulations of the model $\bm{F}_{\bm{\theta}}$.  

To demonstrate the effectiveness of our approach we typically use synthetic (model-generated) data, where we define 
$S_\text{data} = \{ \bm{F}_{\bm{\theta}_0}(\bm{x}_i) : i=1,\ldots,N_\text{data}\}$, where  $\bm{x}_i$ are samples from the initial data distribution and $\bm{\theta}_0$ is given but treated as unknown. In real biological or chemical applications, patterns $S_\text{data}$ come from experimental observations and the unknown $\bm{\theta}_0$ has to be estimated solely from this data. 
The CIL approach allows us to define a stochastic cost function 
that quantifies the discrepancy between the model $\bm{F}_{\bm{\theta}}$ and the pattern data $S_\text{data}$. The actual implementation of the cost function depends on the amount of data $N_\text{data}$.


\subsubsection{Large data sets: Basic Correlation Integral Likelihood}

\begin{figure}
    \centering
    \includegraphics[width=0.99\textwidth]{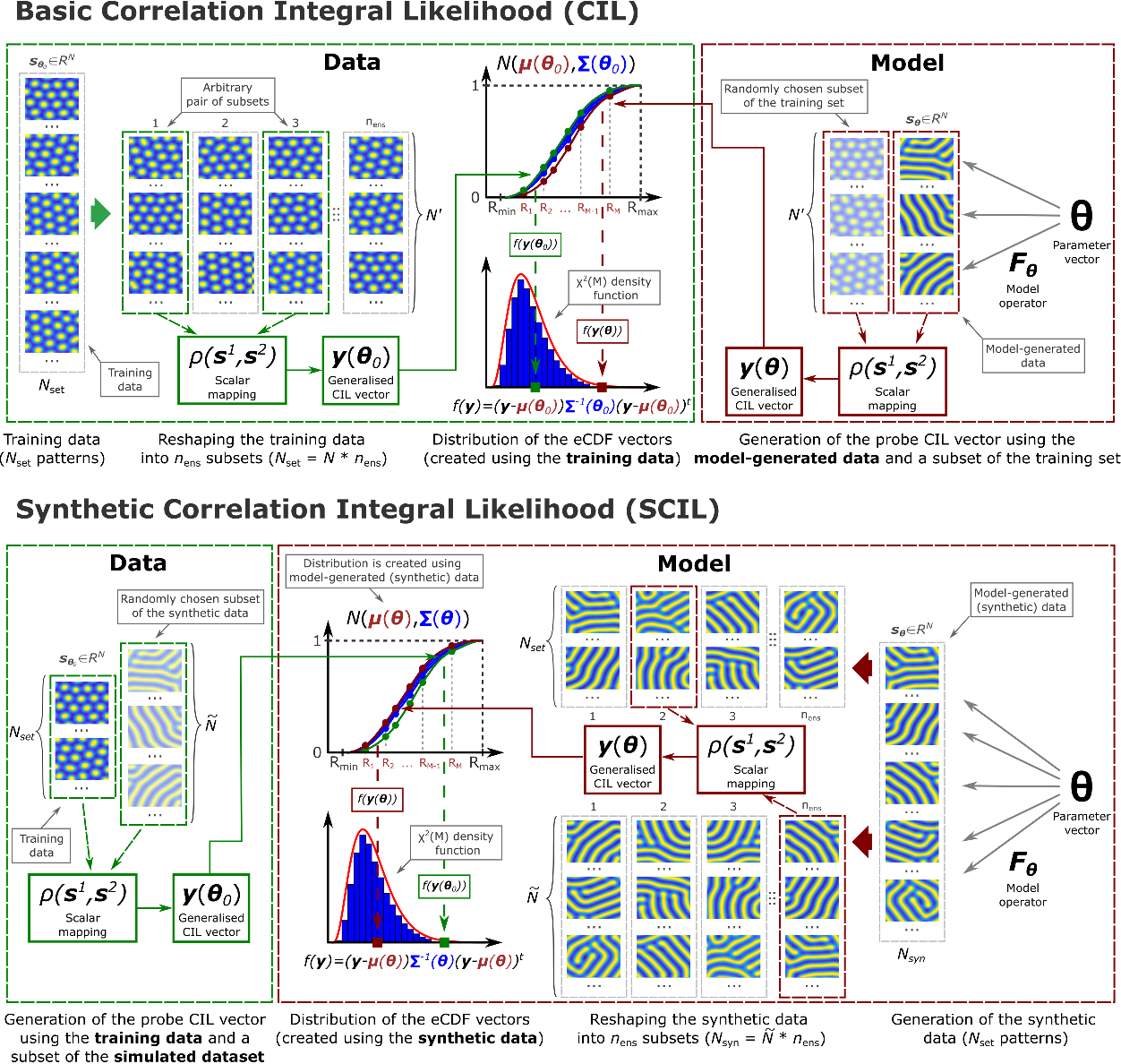}
    \caption{Statistical comparison between a pattern formation model and pattern data using the CIL approach. In the case of basic CIL (top part), the multidimensional Gaussian distribution of $\bm{y}(\bm{\theta}_0)$ is derived from the experimental data. When a new parameter vector $\bm{\theta}$ is proposed, one realisation of $\bm{y}(\bm{\theta})$ is created using model-generated (synthetic) data and compared with the previously derived Gaussian distribution. In the case of SCIL (bottom part) the roles of synthetic and experimental data are reversed: the Gaussian distribution of $\bm{y}(\bm{\theta})$ is created for each proposed parameter vector $\bm{\theta}$ and next compared with a single realisation of $\bm{y}(\bm{\theta}_0)$ derived from the experimental data. }
    \label{fig:cil_schemes}
\end{figure}

If the size  $N_\text{data}$ of the training set is sufficiently large, we can derive enough independent observations of $D_{\bm{\theta}}$ from $S_\text{data}$ to obtain stable estimates for $\bm{\mu}(\bm{\theta}_0)$ and $\bm{\Sigma}(\bm{\theta}_0)$. Producing completely independent observations of distances from unique data sample pairs would typically require a prohibitively large amount of data. Thus, we employ a numerical approximation that allows us to reuse the data samples (Figure~\ref{fig:cil_schemes}, top). We subdivide the training set $S_\text{data}$ into $n_\text{ens}$ subsets, denoted $S^k$, $k=1,\ldots,n_\text{ens}$, each consisting of $N'$ samples of patterns, such that $N_\text{data} = n_\text{ens} \times N'$. Next, computing the eCDF vector of the $N' \times N'$ distances between two sets of $N'$ patterns, $S^k$ and $S^l$, for all possible values of $k \ne l$ provides $\left({n_\text{ens} \atop 2}\right)$ realisations of the correlation integral vector $\bm{y}^{k,l}_{0}$, 
which we use to obtain estimates $\bm{\mu}_0$ and $\bm{\Sigma}_0$ for the parameters of $\bm{y}(\bm{\theta}_0)$. 

\begin{figure}[t]
    \centering
    \includegraphics[width=0.99\textwidth]{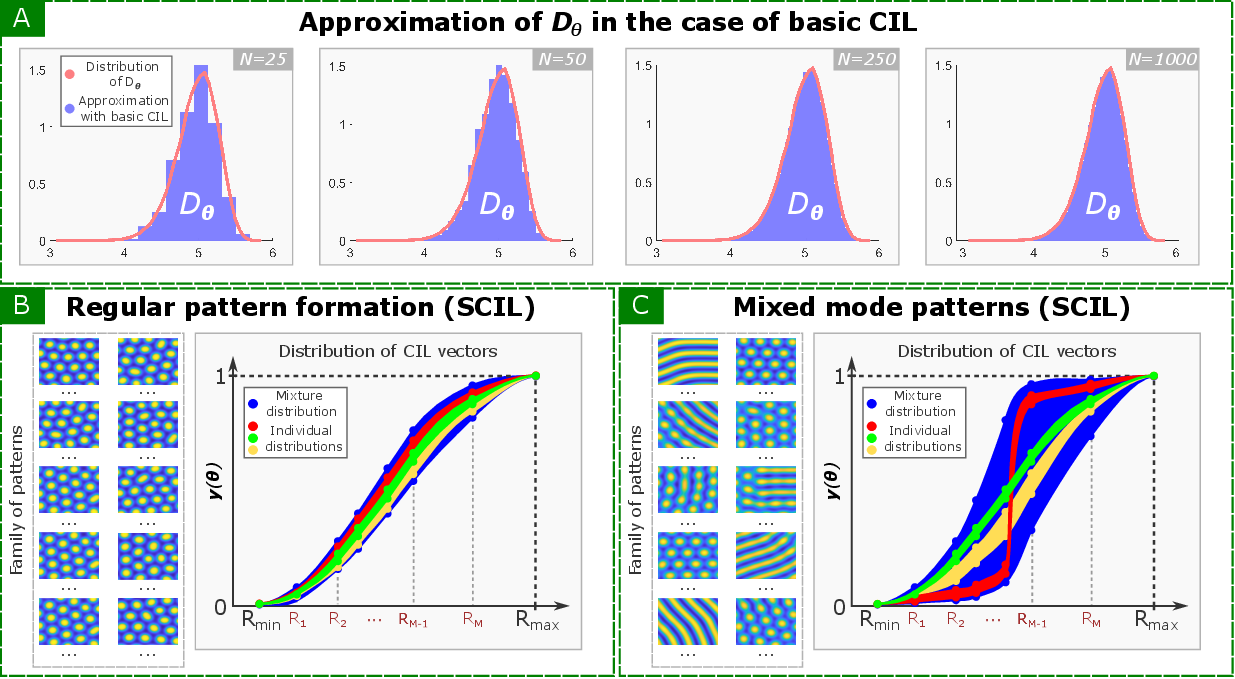}
    \caption{   
    [I]: the comparison of the distribution of distances obtained by basic CIL for different subset size $N$ with the distribution of $D_{\bm{\theta}}$ estimated by computing $10^6$ independent realisations of the random variable. [II]: the distribution of eCDF vectors in the case of SCIL and regular pattern formation. [III]: the distribution of eCDF vectors in the case of SCIL and mixed mode patterns.}
    \label{fig:cil_behaviour}
\end{figure}

In this context, the data is not truly i.i.d.; however, in numerical experiments, this approximation closely follows the density of $D_{\bm{\theta}}$ even for relatively small values of $N'$ (see Figure~\ref{fig:cil_behaviour}, part I). Moreover, according to the theorems of U-statistics, Gaussianity still holds for weakly dependent observations under mild conditions \cite{Borovkova2001,Neumeyer2004}. These conditions, however, are challenging to verify theoretically. Therefore, in practice, we numerically test for Gaussianity of the ensemble of vectors $\{\bm{y}^{k,l}_{0}\}$ using the $\chi^2$-test (or a scalar normality test for each component of the vector). This procedure, combined with bootstrapping, allows for the production of stable estimates for the mean and covariance when $N' \geq 50$. Below this threshold the procedure becomes unstable. 

Once the estimates are obtained, the cost function for parameter estimation can be defined as the negative log-likelihood function
\begin{equation}
	f(\bm{\theta}) = \big(\bm{y}(\bm{\theta})-\bm{\mu}_0 \big)^\top \bm{\Sigma}_0^{-1} \big(\bm{y}(\bm{\theta})-\bm{\mu}_0\big),
	\label{cil:CF}
\end{equation} 
where the correlation integral vector $\bm{y}(\bm{\theta})$ at an arbitrary parameter value $\bm{\theta}$ is evaluated by using $N'$ simulated patterns obtained from $\bm{F}_{\bm{\theta}}$ and one randomly chosen subset $S^k \subset S_\text{data}$. Note that crucially the estimates for $\bm{\mu}_0$ and $\bm{\Sigma}_0$ are constructed off-line from a \enquote{sufficiently large} data set, while the evaluation of $f(\bm{\theta})$ during the parameter estimation only requires a \enquote{small} number $N'$ of model evaluations.

\begin{figure}[t]
    \centering
    \includegraphics[width=0.99\textwidth]{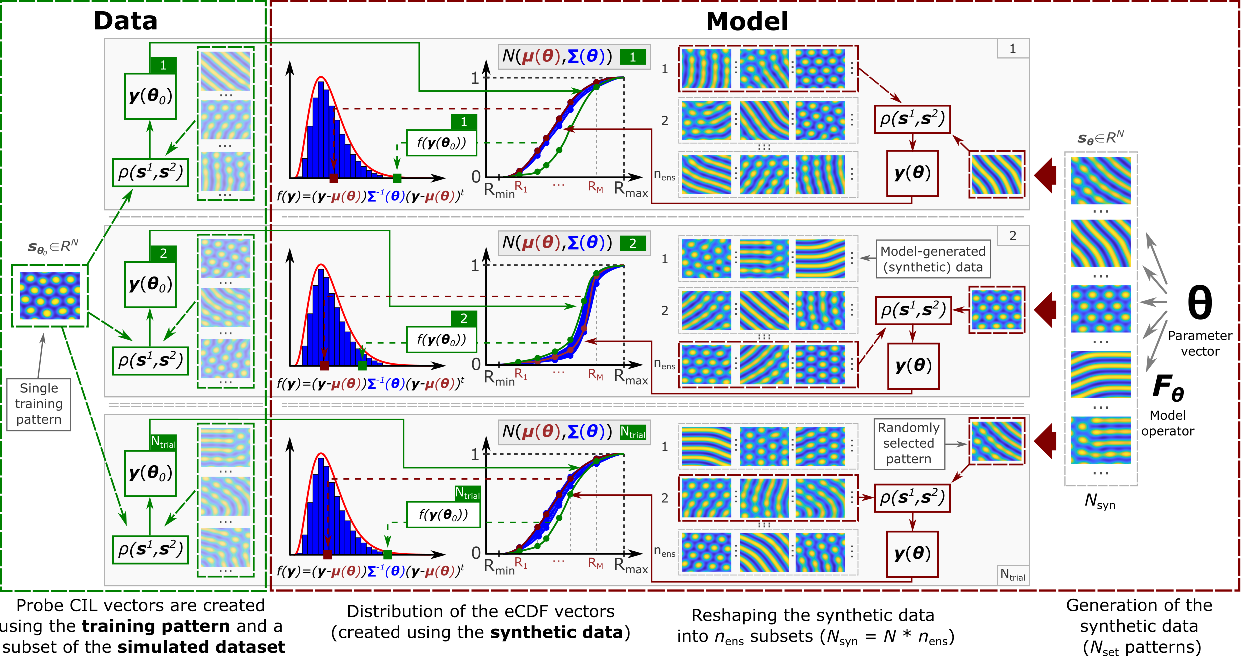}
    \caption{Statistical comparison between a pattern formation model and data using the mixed mode CIL approach. Here, an ensemble of Gaussian distributions $\bm{y}(\bm{\theta})$ defined by fixing small numbers of patterns from the family $\bm{S}_{\bm{\theta}}$ is created for a proposed parameter vector $\bm{\theta}$. Each of these distributions is then compared with a single realisation of $\bm{y}(\bm{\theta}_0)$ derived from the experimental data and the results are processed by summary statistics to produce a single scalar output.}
\end{figure}


    

\subsubsection{Limited data sets: Synthetic Correlation Integral Likelihood (SCIL)}
\label{sec:scil}

When the resampling techniques are insufficient to estimate $\bm{\mu}_0$ and $\bm{\Sigma}_0$ accurately enough from the training set $S_\text{data}$, we can reverse the roles of experimental data and simulations \cite{Kazarnikov2023}. This is of course significantly more expensive during parameter estimation than the original CIL approach. 

To define a cost function $f(\bm{\theta})$ for any given parameter value $\bm{\theta}$, we approximate the CIL at $\bm{\theta}$ using a sufficiently large set $S_\text{syn} = \{\bm{s}^i_{\bm{\theta}} : i=1,\ldots,N_\text{syn}\}$ of patterns, computed from $N_\text{syn}$ model runs, and then compare it to the training data set $S_\text{data}$, as illustrated in Figure~\ref{fig:cil_schemes} (bottom part). This idea is a special case of a \emph{Bayesian synthetic likelihood (BSL) method} \cite{Price2018,Wood2010}.

Since $N_\text{data}$, the number of experimental patterns in $S_\text{data}$ is assumed to be small, we now subdivide the synthetic patterns in $S_\text{syn}$  
into $2 \times n_\text{ens}$ subsets $S^{k,1}_{\bm{\theta}}$ and $S^{k,2}_{\bm{\theta}}$, $k=1,\ldots,n_\text{ens}$, with $N_\text{data}$ and $\widetilde{N}$ patterns each, respectively, such that $(N_\text{data} + \widetilde{N}) \times n_\text{ens} \le N_\text{syn}$. 
By picking $S^{k,1}_{\bm{\theta}}$ and $S^{l,2}_{\bm{\theta}}$, for any $k,l=1,\ldots,n_\text{ens}$, and computing distances between the patterns in the two sets, we can now obtain $n_\text{ens}^2$ realisations  $\widetilde{\bm{y}}^{k,l}_{\bm{\theta}}$ of the synthetic correlation integral vector $\widetilde{\bm{y}}(\bm{\theta})$ at $\bm{\theta}$, which are then used to estimate the parameters $\widetilde{\bm{\mu}}_{\bm{\theta}}$ and $\widetilde{\bm{\Sigma}}_{\bm{\theta}}$ of $\widetilde{\bm{y}}(\bm{\theta})$ (see \cite{Kazarnikov2023} for more details).

We then define the cost function to be 
\begin{equation}
	f(\bm{\theta}) = \big(\bm{y}(\bm{\theta})-\widetilde{\bm{\mu}}_{\bm{\theta}} \big)^\top \widetilde{\bm{\Sigma}}_{\bm{\theta}}^{-1} \big(\bm{y}(\bm{\theta})-\widetilde{\bm{\mu}}_{\bm{\theta}}\big),
	\label{cil:CF1}
\end{equation} 
where $\bm{y}(\bm{\theta})$ is evaluated using the entire training set $S_\text{data}$ and one randomly chosen subset $S^{k,2}_{\bm{\theta}}$ of $\widetilde{N}$ synthetic patterns. Here, unlike the basic CIL, bin values of the eCDF vectors depend on parameter values $\bm{\theta}$ and are automatically estimated during the evaluation of $f(\bm{\theta})$, and here we need $N_\text{syn}$ model runs per cost function evaluation during parameter estimation.




\subsubsection{New approach: Mixed mode Correlation Integral Likelihood}

In the chemical system described in Section \ref{sec:chemical_experiment}, mixed mode patterns appear under 
fixed experimental conditions: One can observe either a pure spot or stripe pattern, or a combination of both behaviours in the same experiment. Provided sufficiently many realisations of experimental patterns are available,  parameter identification can be performed using the basic CIL approach. However, the identification of parameters is more challenging when just a few or a single observation 
is available, since the variability of pattern data is high and might be qualitatively different in the simulations.
In particular, pattern types that are similar to patterns observed in the data may appear in the set of simulated patterns or may not appear at all.  

We deal with this situation by adapting the SCIL approach, where instead of creating a single synthetic likelihood for each $\bm{\theta}$, we create an ensemble of likelihoods.

We generate again a sufficiently large set $S_\text{syn}=\{\bm{s}^i_{\bm{\theta}} : i=1,\ldots,N_\text{syn}\}$ of patterns, computed from $N_\text{syn}$ model runs with parameter $\bm{\theta}$. These patterns are subdivided into $2 \times n_\text{ens}$ subsets $S^{k,1}_{\bm{\theta}}$ and $S^{k,2}_{\bm{\theta}}$, $k=1,\ldots,n_\text{ens}$, each containing $N_\text{data}$ and $\widetilde{N}$ patterns, respectively, with $(N_\text{data} + \widetilde{N}) \times n_\text{ens} \le N_\text{syn}$. We compute the realisations $\widetilde{\bm{y}}^{k,l}_{\bm{\theta}}$, 
$k,l = 1,\ldots,n_\text{ens}$, of the synthetic correlation integral vector at $\bm{\theta}$ as for the SCIL approach in Section \ref{sec:scil}. However, we now proceed slightly differently and define multiple cost functions that will allow us to compute an ensemble of realisations of synthetic likelihoods.
A basic procedure for evaluating the cost function is summarised in Algorithm~\ref{alg:cil_one}. 

\begin{algorithm}[t]
	\SetKwInOut{Input}{Input}\SetKwInOut{Output}{Output}
	\SetKwInOut{Data}{Data}
	\SetKwInOut{Result}{Result}
	\Data{$S_\text{data} \,\ldots\ $ a set of patterns $\bm{s}^1,\ldots,\bm{s}^{N_\text{data}}$ with unknown model 
    parameter $\bm{\theta}_0$}
	\Input{$\rho_\alpha  \quad\ \ldots\ $ a family of distances with $\alpha=1,\ldots,N_\text{dist}$}
	\Input{$\bm{\theta} \quad\ \ \, \ldots$ \ the parameter value where the cost function should be evaluated}
	\Input{$N_\text{syn} \ \, \ldots$ \ the number of model-generated (synthetic) patterns at $\bm{\theta}$}
        \Input{$G \quad \ \ \, \ldots$ \ the summary statistics (e.g., $G= \min$ as in \eqref{eq:minfk})}
	\Output{$f(\bm{\theta}) \ \; \ldots$ \ the value of the  cost function at $\bm{\theta}$}	
	\Begin
	{
		1. Simulate $N_\text{syn}$ patterns $\bm{s}^i_{\bm{\theta}}$, $i=1,\ldots,N_\text{syn}$, with parameter $\bm{\theta}$
		
		2. Divide these patterns into $2 \times n_\text{ens}$ subsets $S_{\bm{\theta}}^{k,1}$ and $S_{\bm{\theta}}^{k,2}$, $k=1,\ldots,n_\text{ens}$, with $N_\text{data}$\\ \quad \ and $\widetilde{N}$ patterns, respectively

        3. \For{$k=1,\ldots,n_\mathrm{ens}$}
        {
            3.1 \For{$l=1,\ldots,n_\mathrm{ens}$}
		{
		3.1.1 Initialise the correlation integral vector $\bm{y}_{\bm{\theta}}^{k,l}$ for $S_{\bm{\theta}}^{k,1}$ and $S_{\bm{\theta}}^{l,2}$ to be empty\!\!\! 
  
		3.1.2 \For{$\alpha=1,\ldots,N_\mathrm{dist}$}
		    {
			* Compute the distances $\rho_\alpha(\bm{s}^1,\bm{s}^2)$ between all patterns $\bm{s}^1 \in S_{\bm{\theta}}^{k,1}$ and $\bm{s}^2 \in S_{\bm{\theta}}^{l,2}$\!\!\!\!
		
			* Compute part of correlation integral vector $\bm{y}^{k,l}_{\bm{\theta},\alpha}$ for $\rho_\alpha$
			
			* Concatenate the current vector $\bm{y}_{\bm{\theta}}^{k,l}$ and $\bm{y}^{k,l}_{\bm{\theta},\alpha}$
		    } 
            } 
            3.2 Using the samples $\bm{y}_{\bm{\theta}}^{k,l}$ computed in Step 3.1, estimate $\bm{\mu}_{\bm{\theta}}^k$ and $\bm{\Sigma}_{\bm{\theta}}^k$ of the (multi-feature) CIL vector 
            corresponding to the fixed set of patterns $S_{\bm{\theta}}^{k,1}$ 
            at $\bm{\theta}$ 
            
            3.3 Randomly select one of the subsets $S_{\bm{\theta}}^{l,2}$ of $\widetilde{N}$ patterns from $\bm{s}^i_{\bm{\theta}}$, $i=1,\ldots,N_\text{syn}$
            
            3.4 Using the computed estimates $\bm{\mu}_{\bm{\theta}}^k$, $\bm{\Sigma}_{\bm{\theta}}^k$ for the fixed set of patterns $S_{\bm{\theta}}^{k,1}$ and the set $S_{\bm{\theta}}^{l,2}$, evaluate the cost function $f^k(\bm{\theta})$ in \eqref{cil:CF2} with
            $\bm{y}(\bm{\theta})$ as defined after~\eqref{cil:CF2}
            } 
		
        4. Process the scalar values $f^k(\bm{\theta})$ computed in Step 3 with the summary statistics $G$ and evaluate the cost function at $\bm{\theta}$ as follows: $f(\bm{\theta}) = G(\{ f_1(\bm{\theta}), \ldots, f_{N_\text{trial}}(\bm{\theta})\})$.
		
	}
	\caption{Construction and evaluation of  $f(\bm{\theta})$ for mixed mode SCIL at a single parameter value $\bm{\theta}$}
	\label{alg:cil_one}
\end{algorithm}

We fix $k \in \{1,\dots,n_\text{ens}\}$ and consider all realisations $\widetilde{\bm{y}}^{k,l}_{\bm{\theta}}$, $l = 1,\ldots,n_\text{ens}$, computed by comparing the $N_\text{data}$ fixed patterns in $S^{k,1}_{\bm{\theta}}$ with the patterns in the sets $S^{l,2}_{\bm{\theta}}$. We estimate the mean $\widetilde{\bm{\mu}}^k_{\bm{\theta}}$ and covariance $\widetilde{\bm{\Sigma}}^k_{\bm{\theta}}$ of these $n_\text{ens}$ realisations and define for each choice of the (synthetic) $N_\text{data}$ sample patterns $S^{k,1}_{\bm{\theta}}$ a separate negative log-likelihood
\begin{equation}
	f^k(\bm{\theta}) = \big(\widetilde{\bm{y}}(\bm{\theta})-\bm{\mu}_{\bm{\theta}}^k \big)^\top \left(\bm{\Sigma}_{\bm{\theta}}^k\right)^{-1} \big(\widetilde{\bm{y}}(\bm{\theta})-\bm{\mu}_{\bm{\theta}}^k\big),
	\label{cil:CF2}
\end{equation} 
where $\widetilde{\bm{y}}(\bm{\theta})$ is computed by comparing the experimental patterns in the training set $S_\text{data}$ with the synthetic patterns in one randomly chosen subset $S^{l,2}_{\bm{\theta}}$, $l = 1,\ldots,n_\text{ens}$.

As the $N_\text{data}$ patterns in $S^{k,1}_{\bm{\theta}}$ that are used to compute $f^k(\bm{\theta})$ are kept fixed here, the associated likelihood will depend on them. This situation is illustrated in Figure~\ref{fig:cil_behaviour} (parts II and III), showing how the pattern distribution $\bm{S}_{\bm{\theta}}$ gives rise to an ensemble of Gaussian distributions, defined by different choices of the $N_\text{data}$ synthetic patterns used to compute them. In the regular SCIL setting, i.e., when all the patterns are qualitatively similar, the difference between the individual distributions is not too pronounced. 
In the mixed mode situation, on the other hand, we obtain widely different sets of SCIL vectors.

The final likelihood to be used in our parameter estimation for mixed mode patterns is now derived using 
\begin{equation}
\label{eq:minfk}
f(\bm{\theta}) = \min\limits_{1\leq k \leq n_\text{ens}}f^k(\bm{\theta}),
\end{equation}
i.e., a parameter value $\bm{\theta}$ should be \enquote{accepted}, if the simulations with that value may produce patterns that are close to the patterns in $S_\text{data}$. With this selection, we arrive at the extreme value distribution \cite{Fisher1928} to be used as the cost function. Other  applications may call for other choices here, e.g., giving a high likelihood to a parameter value only if sufficiently many simulated patterns are close to the patterns in the training set.
%
%
%

Depending on the required and available computational resources, the approach may be implemented in different ways. Since here only $n_\text{ens}$ realisations of $\widetilde{\bm{y}}^{k,l}_{\bm{\theta}}$ are used to estimate the parameters of the cost function $f^k(\bm{\theta})$ in \eqref{cil:CF2}, as opposed to $n_\text{ens}^2$ realisations in the original SCIL approach, the total number $N_\text{syn}$ of model runs needs to be significantly larger. In addition, bootstrapping may be applied to reduce the required number of model simulations, significantly decreasing computational cost without losing accuracy (see Algorithm~\ref{alg:CIL3bst}).


While mixed mode SCIL is a robust general solution for the case when $N_\text{data} < \widetilde{N}$, sometimes a simplified scheme of SCIL can be used as a less computationally expensive alternative.  If the bias created by fixed patterns is relatively small, the SCIL approach gives a sufficiently accurate approximation of the distribution of $D_{\bm{\theta}}$ (see Figure~\ref{fig:cil_behaviour} (part II)). This was the situation studied in our previous work \cite{Kazarnikov2023}, where we focused on synthetic data and regular pattern formation, where each set of model parameters produced a single type of pattern (e.g., stripes, labyrinths, or hexagons). 

\begin{remark}
One can show that the realisations $\widetilde{\bm{y}}^{k,l}_{\bm{\theta}}$ of the eCDF vector in the mixed-mode SCIL approach (for fixed $k$) still follow a Gaussian distribution. 

Let us draw $N_\text{data}$ patterns $\bm{s}^1,\ldots,\bm{s}^{N_\text{data}}$ from the random family $\bm{S}_{\bm{\theta}}$ and consider for each fixed pattern $\bm{s}^j$ the i.i.d. scalar random variables $D_{\bm{\theta}}^{j,i} = \rho(\bm{s}^j, \bm{S}^i_{\bm{\theta}})$, where $\bm{S}^i_{\bm{\theta}}$, $i=1,\ldots,n$, are i.i.d. copies of $S_{\bm{\theta}}$. 
For fixed bin values $R_1,\ldots,R_M$, we now define the random vector $\widetilde{\bm{Y}}$ with components
\[
\widetilde{Y}_k = \dfrac{1}{N_\text{data}}
\sum\limits_{j=1}^{N_\text{data}}
\left(
\dfrac{1}{n}
\sum\limits_{i=1}^n
\#(D_{\bm{\theta}}^{j,i} < R_k)
\right),
\quad
k=1,\ldots,M.
\]
Using the CLT and the fact that a sum of Gaussians is again Gaussian, 
we can show that asymptotically, for $n \to \infty$, $\widetilde{\bm{Y}}$ follows a multivariate Gaussian distribution. 
Here, the random distances $D_{\bm{\theta}}^{j,i}$ are evaluated by using randomly drawn patterns from $S_{\bm{\theta}}$ in the second argument. However, the mean and variance of this Gaussian distribution naturally depend on the $N_\text{data}$ fixed patterns $\bm{s}^1,\ldots,\bm{s}^{N_\text{data}}$ in the first argument.
\end{remark}


\begin{algorithm}[t]
	\SetKwInOut{Input}{Input}\SetKwInOut{Output}{Output}
	\SetKwInOut{Data}{Data}
	\SetKwInOut{Result}{Result}
	\Data{$S_\text{data} \,\ldots\ $ a set of patterns $\bm{s}^1,\ldots,\bm{s}^{N_\text{data}}$ with unknown model 
    parameter $\bm{\theta}_0$}
	\Input{$\rho_\alpha  \quad\ \ldots\ $ a family of distances with $\alpha=1,\ldots,N_\text{dist}$}
	\Input{$\bm{\theta} \quad\ \ \, \ldots$ \ the parameter value where the cost function should be evaluated}
	\Input{$N_\text{syn} \ \, \ldots$ \ the number of model-generated (synthetic) patterns at $\bm{\theta}$}
    \Input{$N_\text{trial} \ \ldots$ \ the number of individual distributions to sample}
	\Input{$n_\text{CIL} \ \ \ldots$ \ \ Number of generalised correlation integral vectors to estimate MCIL}
    \Input{$G \quad \ \ \, \ldots$ \ the summary statistics (e.g., $G= \min$ as in \eqref{eq:minfk})}
	\Output{$f(\bm{\theta}) \ \; \ldots$ \ the value of the  cost function at $\bm{\theta}$}        
	\Begin
	{
		1. Simulate $N_\text{syn}$ patterns $S_\text{syn} := \{\bm{s}^i_{\bm{\theta}}:i=1,\ldots,N_\text{syn}\}$, with parameter $\bm{\theta}$

            2. \For{$k=1,\ldots,N_\text{trial}$}
            {
                2.1 Construct a subset $S^1$ by randomly selecting  $N_\text{data}$ patterns from $S_\text{syn}$ (without replacement)      
                
                2.2 \For{$l=1,\ldots,n_\mathrm{CIL}$}
		      {		      	
		          2.2.1 Construct the subset $S^2$ by randomly selecting with replacement $\widetilde{N} = N_\text{syn} - N_\text{data}$ patterns from the remaining patterns in $S_\text{syn}$
		
		          2.2.2 Initialise the correlation integral vector $\bm{y}_{\bm{\theta}}^{l}$ to be the empty vector

		          2.2.3 \For{$\alpha=1,\ldots,N_\mathrm{dist}$}
		    {
			* Compute the distances $\rho_\alpha(\bm{s}^1,\bm{s}^2)$ between all patterns $\bm{s}^1 \in S^{1}$ and $\bm{s}^2 \in S^{2}$
			
			* Compute the part of the correlation integral vector $\bm{y}^{l}_{\bm{\theta},\alpha}$ for $\rho_\alpha$
			
			* Concatenate the current vector $\bm{y}_{\bm{\theta}}^{l}$ and $\bm{y}^{l}_{\bm{\theta},\alpha}$
		    } 
                } 

                2.3 Using the samples $\bm{y}_{\bm{\theta}}^{l}$, $l=1,\ldots n_\mathrm{CIL}$ computed in Step 2.2, estimate $\bm{\mu}_{\bm{\theta}}^k$ and $\bm{\Sigma}_{\bm{\theta}}^k$ of the (multi-feature) CIL vector 
            corresponding to the fixed set of patterns $S^{1}$ 
            at $\bm{\theta}$ 
            
                2.4 Randomly select a subset $S^{2}$ of $\widetilde{N}$ patterns from $S_\text{syn}$ (without replacement) 
                
                2.5 Using the computed estimates $\bm{\mu}_{\bm{\theta}}^k$, $\bm{\Sigma}_{\bm{\theta}}^k$ for the fixed set of patterns $S^{1}$ and set $S^{2}$ from Step 2.4, evaluate the cost function $f^k(\bm{\theta})$ in \eqref{cil:CF2} with
            $\bm{y}(\bm{\theta})$ defined after~\eqref{cil:CF2}\!\!\!
            } 
		
            3. Process the scalar values $f^k(\bm{\theta})$ computed in Step 3 with the summary statistics $G$ and evaluate the cost function at $\bm{\theta}$ as follows: $f(\bm{\theta}) = G(\{ f_1(\bm{\theta}), \ldots, f_{N_\text{trial}}(\bm{\theta})\})$.	
	}
	\caption{Using bootstrapping to evaluate $f(\bm{\theta})$ for mixed mode SCIL at a single parameter value $\bm{\theta}$\!\!\!}
	\label{alg:CIL3bst}
\end{algorithm}

\section{Applications to the CIMA reaction \label{sec:parameter_identification}}


In this section, we illustrate the efficacy of our statistical approach in identifying the parameters of the Lengyel-Epstein model \eqref{eq:LEmodel} for the CIMA reaction. As data we use
two experimental patterns from the previously published chemical experiments in~\cite{Rudovics1996}, which are shown in Figure~\ref{fig:chemical_experiment}.

\begin{figure}
    \centering
    \includegraphics[width=0.99\textwidth]{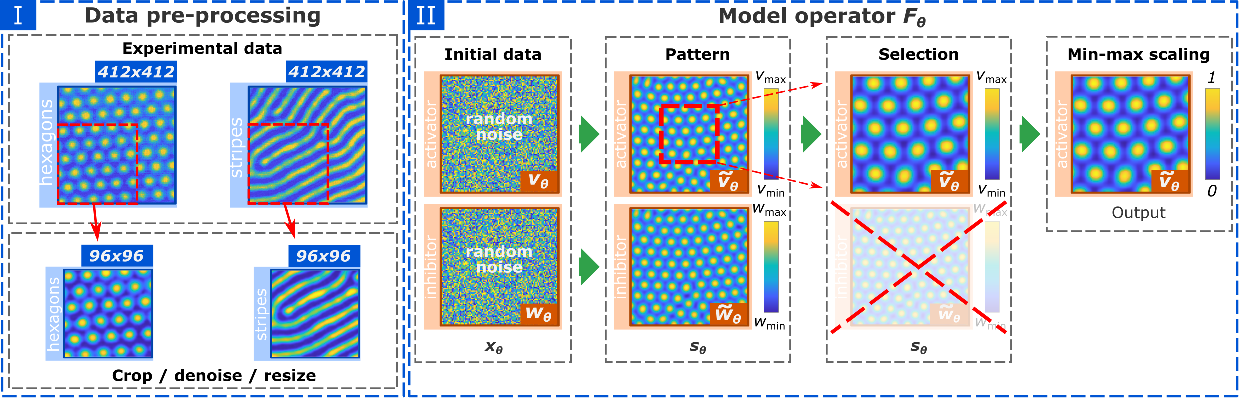}
    \caption{Main steps in setting up the CIL approach to fit Lengyel-Epstein to experimental pattern data: Pre-processing digital images of experimental patterns [I]; setup of the model operator $\bm{F}_{\bm{\theta}}$ [II].}
    \label{fig:model}
\end{figure}

We begin by pre-processing the image data.  Both pattern images have a resolution of 412~x~412 pixels and contain approximately 10 spatial \enquote{wavelengths}.  Since a high number of wavelengths requires more pixels to represent the data accurately, we may need to use high spatial resolution in model simulations to produce comparable synthetic patterns.  To optimise the computational complexity of the parameter identification, we crop each image by extracting a square patch of 275~x~275 pixels.  Next, we remove the measurement noise from the data by smoothing the pixel values with a Gaussian filter and rescale the resulting images to a resolution of 96 by 96 pixels. This procedure is illustrated in Figure~\ref{fig:model}, left.

\subsection{Setup of the numerical model}

To apply the statistical approach discussed in the previous section, we need to define the model operator $\bm{F}_{\bm{\theta}}$, the parameter vector $\bm{\theta}$, and the distribution of initial conditions $\bm{X}(\omega)$ in terms of the deterministic Lengyel-Epstein model \eqref{eq:LEmodel}. We begin by discretising the PDE system \eqref{eq:LEmodel} using the Method of Lines. To eliminate the influence of boundary conditions on the numerically computed patterns, we consider the system within a larger computational domain $\widetilde{\Omega} = (0, \widetilde{L})^2$ with $\widetilde{L} = 2L$ to ensure that $\Omega_0 \subset \widetilde{\Omega}$. We discretise $\widetilde{\Omega}$ using a uniform grid with a fixed step size $h=\widetilde{L}/(M_\text{dim}-1)$, where $M_\text{dim}\in \mathbb{N}$, leading to a finite set of grid points: 
\[
\bm{\xi}_{i,j}=((i-1)h,(j-1)h), \quad i,j=1,\dots,M_\text{dim}.
\]

To reduce the reaction-diffusion system \eqref{eq:LEmodel} to a finite set of $2M_\text{dim}^2$ ordinary differential equations (ODEs), we define for each grid point the time-dependent functions
\[
v_{i,j}(t) = v(\bm{\xi}_{i,j},t),
\ \ 
w_{i,j}(t) = w(\bm{\xi}_{i,j},t),
\qquad
i,j=1,\ldots,M_\text{dim},
\]
and discretise the Laplace operator by the five-point stencil \cite{Grossman2007}
\[
\Delta v\approx\frac{v_{i+1,j}+v_{i-1,j}+v_{i,j+1}+v_{i,j-1}-4v_{i,j}}{h^{2}}=\nabla_{h}^{2}v_{i,j}.
\]
Substituting the approximations leads to the system:
\begin{eqnarray}
	{\textstyle \frac{\text{d}}{\text{d}t}} v_{i,j}(t) & = & \dfrac{1}{\sigma}\left(a - v_{i,j}(t) - 4\dfrac{v_{i,j}(t)w_{i,j}(t)}{1+v_{i,j}(t)^2} + \nabla_{h}^{2} v_{i,j}(t) \right), \nonumber \\[1ex] 
	{\textstyle \frac{\text{d}}{\text{d}t}} w_{i,j}(t) & = & b \left( v_{i,j}(t) - \dfrac{v_{i,j}(t)w_{i,j}(t)}{1+v_{i,j}(t)^2} \right) + d \nabla_{h}^{2}w_{i,j}(t), \quad \label{eq:MOLSystem}
	i,j=1,2,\dots,M_{dim}.
\end{eqnarray}
The Neumann boundary conditions are taken into account by using a one-sided first-order difference scheme \cite{Grossman2007}, which leads to the following conditions for the \enquote{ghost} values $\bm{\xi}_{0,j}$, $\bm{\xi}_{N+1,j}$, $\bm{\xi}_{i,0}$, $\bm{\xi}_{i,N+1}$:
\begin{eqnarray*}
	& v_{0,j}(t) \equiv v_{1,j}(t),
	\,
	v_{M_{dim}+1,j}(t) \equiv v_{M_{dim},j}(t), 
        \,
        v_{i,0}(t) \equiv v_{i,1}(t),
	\,
	v_{i,M_{dim}+1}(t) \equiv v_{i,M_{dim}}(t), \\
	& w_{0,j}(t) \equiv w_{1,j}(t),
	\,
	w_{M_{dim}+1,j}(t) \equiv w_{M_{dim},j}(t),
	\,
	w_{i,0}(t) \equiv w_{i,1}(t),
	\,
	w_{i,M_{dim}+1}(t) \equiv w_{i,M_{dim}}(t).
\end{eqnarray*}

The model parameters in the discretised system \eqref{eq:MOLSystem} are grouped into the vector $\bm{\theta} = (L, a, b, \sigma)$. The parameter $d$, i.e., the ratio of the diffusion coefficients, is fixed to the experimentally known value $d = 1.07$ (see \cite{Rudovics1999} for details). Here, we treat the size of the computational domain $L$ as a free parameter because the model variables are defined in rescaled dimensionless form, and $L$ encodes information about one of the constant reactant concentrations. The distribution of the initial conditions $\bm{X}(\omega)$ depends on $\bm{\theta}$ and it is derived from the conditions of the chemical experiment described in \cite{Rudovics1996}. As the contents of both reservoirs of the chemical reactor are continuously mixed, we prescribe the initial data by the values of the spatially homogeneous steady state of the Lengyel-Epstein model \eqref{eq:LEmodel}, 
$(v_0(\bm{\theta}),w_0(\bm{\theta})) = (a/5, 1 + (a/5)^2)$. This is then perturbed with small concentration fluctuations to initiate chemical pattern formation, which leads to the distribution
 \begin{align}
 \label{eq:initial_data}
 \bm{X}_{\bm{\theta}}(\omega) &= (\bm{V}_{\bm{\theta}}(\omega), \bm{W}_{\bm{\theta}}(\omega)), \quad \text{where}\\
\,
 V_{\bm{\theta}}^{i,j} &\sim v_0(\bm{\theta}) + {U(-\delta, \delta)},
 \ \
 W_{\bm{\theta}}^{i,j} \sim w_0(\bm{\theta}) + {U(-\delta, \delta)},
 \quad
 i,j=1,\ldots,M_\text{dim},\nonumber
 \end{align}
 and $U(-\delta,\delta)$ is uniform random noise with $\delta = 10^{-3}$.
 
 The evaluation of the model operator $\bm{F}_{\bm{\theta}} (\bm{x})$ involves the following steps (Figure~\ref{fig:model}, right):
 \begin{enumerate}
     \item A realisation of randomised initial data $\bm{x}_{\bm{\theta}} \in \mathbb{R}^{2M_\text{dim}^2}$ is generated using formula \eqref{eq:initial_data};
     \item The initial data $\bm{x}_{\bm{\theta}}$ are propagated to the steady-state pattern $\bm{s}_{\bm{\theta}}=(\widetilde{\bm{v}}, \widetilde{\bm{w}})$ by numerically integrating the MOL system \eqref{eq:MOLSystem}. For this purpose, we use the explicit stabilised Runge-Kutta method ROCK2 \cite{Abdulle2001}. The numerical solution is computed over the time interval $[0,T_\text{end}]$ or until the transient behaviour terminates. This is determined by monitoring the $L_2$-norm of the time derivative and stopping the simulation if this value falls below a small tolerance $\varepsilon_\text{conv}$. To reduce numerical fluctuations, the value of the $L_2$-norm is smoothed by computing a moving average over a sufficiently large number of time steps $N_\text{ma}$. If the $L_2$-norm of the time derivative still exceeds $5\varepsilon_\text{conv}$ at $T_\text{end}$, the parameter value is rejected;
     \item A square block $\widetilde{\bm{v}}^c \in \mathbb{R}^{(M_\text{dim}/2)^2}$ of values corresponding to a patch in $\Omega_c \subset \widetilde{\Omega}$ of size $L \times L$ is extracted from the centre of the activator pattern $\widetilde{\bm{v}}$;
     \item The vector $\widetilde{\bm{v}}^c$ is scaled by min-max normalisation to obtain the final, scaled pattern $\bm{s}$ with
     \[s_{i,j}  = (\widetilde{v}^c_{i,j} - \widetilde{v}^c_\text{min}) / (\widetilde{v}^c_\text{max} - \widetilde{v}^c_\text{min}),
     \quad
     \widetilde{v}^c_\text{max} = \max\limits_{1\leq i,j\leq M_\text{dim}/2} \widetilde{v}^c_{i,j}\,,
     \quad
     \widetilde{v}^c_\text{min} = \min\limits_{1\leq i,j\leq M_\text{dim}/2} \widetilde{v}^c_{i,j}\,.
     \]
 \end{enumerate}

\subsection{Application of the mixed mode SCIL approach to the CIMA model}

We define the CIL cost function $f(\bm{\theta})$, using Algorithm~\ref{alg:CIL3bst} with a family of norms suggested in our previous work \cite{Kazarnikov2023} for min-max scaled data, namely $\| \cdot \|_{L^2}$, $\| \cdot \|_{W^{1,2}}$, and $\vvvert \cdot \vvvert_{W^{1,2}}$. 

The mixed mode SCIL approach now consists in first minimising the cost function with respect to the parameter vector $\bm{\theta}$ using Differential Evolution (DE) \cite{Storn1997}. All relevant parameters for the numerical experiment are provided in Table~\ref{tab:experiment_settings}. DE begins with a \enquote{population} of potential candidate parameters and mimics the process of evolution by iteratively improving the population to reach a local optimum of the optimisation problem.  The initial DE parameter populations are drawn from $L \sim U(25,80)$, $a \sim U(5,15)$, $b \sim U(0.05,0.5)$ and $\sigma \sim U(1,100)$. By working with multiple candidates rather than focusing on a single parameter, the algorithm fosters diversity and better explores the solution space. DE involves three key operations: mutation, which introduces new trial vectors by perturbing existing solutions with a scaled difference between randomly selected candidates; crossover, which combines information from the trial vectors and the original population to create offspring, and selection, where DE evaluates the fitness of trial vectors and selects the best ones to form the next generation. We iterate the method until convergence to the local minimum is achieved, and the cost function values $f(\bm{\theta})$ no longer improve.

\begin{table}
    \centering
    \begin{tabular}{|c|c|}
    \hline
         \textbf{Parameter name} & \textbf{Value}  \\\hline
        $M_\text{dim}$ (spatial resolution) & 128 \\
         $T_\text{end}$ (end time for integration interval) & $2 \times 10^6$ \\
         $\varepsilon_\text{conv}$ (convergence threshold) & $2 \times 10^{-7}$
         \\\hline
         $N_\text{syn}$ (number of synthetic patterns at $\bm{\theta}$) & 800 \\
         $N_\text{trial}$ (number of distributions in ensemble) & 100 \\
         $n_\text{CIL}$ (number of samples to estimate CIL) & 1000 \\
         $M$ (dimension of the CIL vector) & $3 \times 12$ \\
         \hline
         DE population number & 39 \\
         crossover probability & $0.9$ \\
         differential weight & $0.8$ \\\hline
    \end{tabular}
    \caption{Parameter values used in the numerical experiments.}
    \label{tab:experiment_settings}
\end{table}

\begin{figure}
    \centering
    \includegraphics[width=0.99\textwidth]{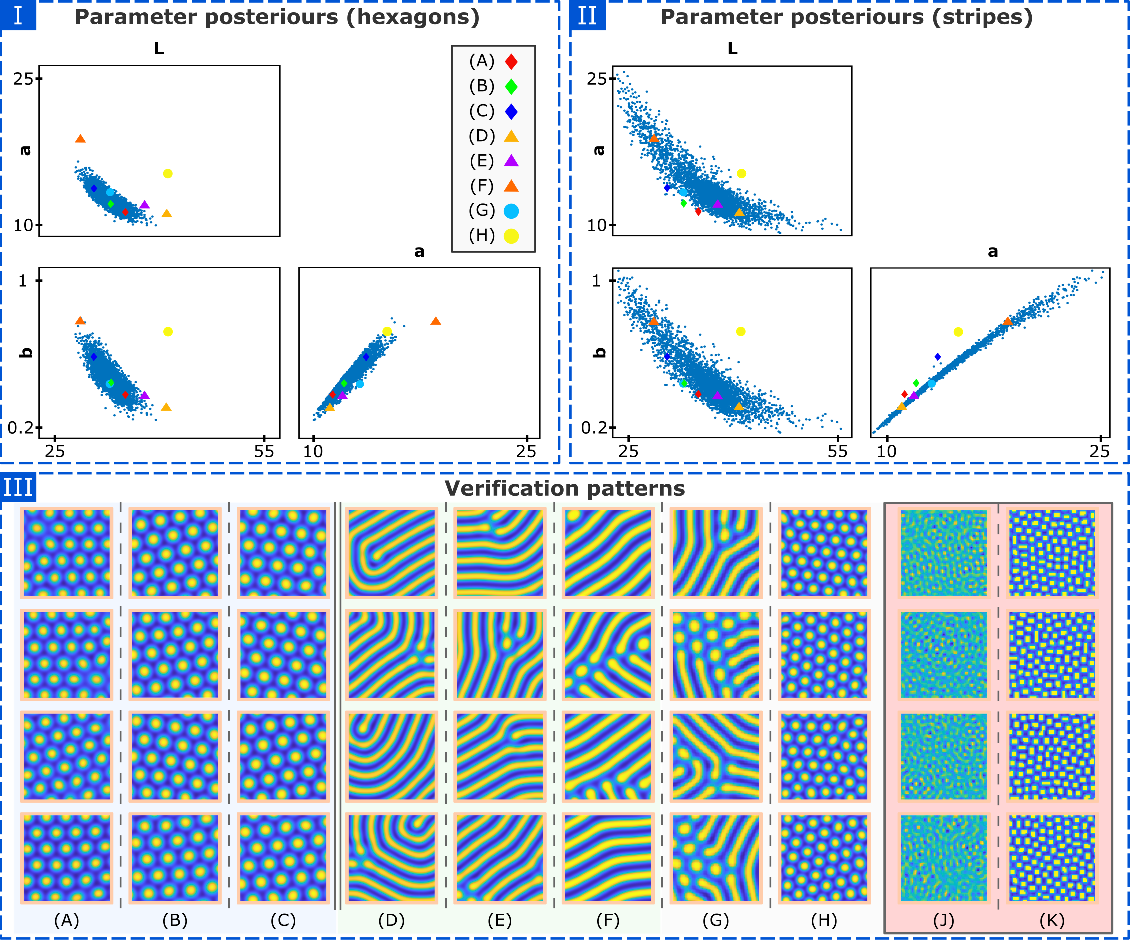}
    \caption{Parameter identification results obtained with the mixed mode SCIL approach: Parameter posterior distributions corresponding to hexagonal and striped experimental patterns [I] and [II] with 8 verification parameter values and respective simulated patterns [III]. The final two simulated patterns (J) and (K) correspond to the model parameters derived from chemical measurements for the hexagonal and striped patterns, respectively. These parameter values are shown in Table~\ref{tab:parameter_estimates}.}
    \label{fig:results}
\end{figure}

Subsequently, we proceed with uncertainty quantification using the Adaptive Metropolis sampling method \cite{Haario2001}, utilising the best candidate of the final DE population as the starting point for the Markov Chain Monte Carlo (MCMC) sampler. 
As our final parameter estimate we use the maximum a posteriori 
(MAP) point of the empiricial distribution. The MAP estimates for the two esperimental patterns are presented in Table~\ref{tab:parameter_estimates} together with 95\% credible intervals, while parameter posteriors are plotted in Figure~\ref{fig:results} together with verification parameter values used to validate the algorithm performance. It can be seen from the validation points that the statistical approach is able to correctly separate different types of model behaviour, showing an accuracy comparable with the \enquote{eyeball norm}. The hexagonal patterns (A)-(C) all lie clearly within the centre of the posterior distribution in Panel [I], while the striped patterns (D)-(F) lie outside the empirical distribution. For the posterior distribution in panel [II] the situation is exactly reversed. The parameter value for the mixed mode patterns in (G), on the other hand, lies at the edge of both posterior distributions. The patterns with slightly smaller hexagons in (H) are outside of both distributions.

The variability with respect to $\sigma$ is significantly higher compared to the other parameters (see Table~\ref{tab:parameter_estimates}). This is explained by the fact that this parameter is proportional to the concentration of the complexing agent $[\text{S}]$, which mainly impacts the transient phase and the time needed for patterns to form. However, our parameter identification approach only uses the final stationary pattern and we do not take into account differences in how long it takes the model to reach the steady state at a particular parameter value.  Since $\sigma$ is not the parameter of interest, we treat it as a nuisance parameter and integrate it out, considering only the posterior of the restricted parameter $\bm{\theta}' = (L,a,b)$ in Figure~\ref{fig:results}. 

\begin{table}
    \centering
\begin{tabular}{|l|c|c|c|c|}
\hline
 & \textbf{L} & \textbf{a} & \textbf{b} & \boldmath$\sigma$ \\
\hline
\multicolumn{5}{|c|}{\textbf{Hexagonal pattern}} \\
\hline
MAP estimates & 33.7 & 12.7 & 0.45 & 385 \\
\hline
95\% Credible Interval & $[31.0, 36.8]$ & $[11.7, 15.0]$ & $[0.37, 0.64]$ & $[131, 4557]$ \\
\hline
Experimental estimates & 57.5 & 23.3 & 1.61 & 151 \\
\hline
\multicolumn{5}{|c|}{\textbf{Striped pattern}} \\
\hline
MAP estimates & 38.1 & 12.6 & 0.39 & 162 \\
\hline
95\% Credible Interval & $[27.9, 46.0]$ & $[11.0, 20.7]$ & $[0.28, 0.85]$ & $[62, 2158]$ \\
\hline
Experimental estimates & 57.5 & 25.9 & 1.61 & 151 \\
\hline
\end{tabular}

    \caption{MAP parameter estimates with 95\% credible intervals obtained from the CIL approach and parameter values derived from chemical measurements.} 
    \label{tab:parameter_estimates}
\end{table}

The numerical integration of the discretised model \eqref{eq:MOLSystem} was performed using the previously developed efficient parallel implementation of the ROCK2 numerical solver. We employed modern Graphical Processing Units (GPUs), using the Nvidia CUDA (Compute Unified Device Architecture) computing platform to execute computations on massively parallel Nvidia GPU devices. The code used for the numerical experiments is freely available at \url{https://github.com/AlexeyKazarnikov/CILNumericalCode}. All GPU-based code has been developed with CUDA 11.0 and compiled with Microsoft Visual C++ 2019 under Microsoft Windows and GCC 9.3.0 under Ubuntu Linux. 



\renewcommand{\topfraction}{.8}
\renewcommand{\floatpagefraction}{.8}

\subsection{Comparison to experimental estimates derived from chemical measurements}

To put our parameter estimates in context, we compare them with experimental estimates derived from chemical measurements reported in \cite{Rudovics1996}. This derivation is comprised of multiple stages, each characterised by a significant degree of uncertainty. Since reaction rates depend on temperature, and the temperature reported in \cite{Rudovics1996} is $5 \pm 0.1 ^\circ C$, we approximate them using estimates for $4 ^\circ C$ from various experimental studies \cite{Rudovics1999} (see Table~\ref{tab:chemical_concentrations}, left). The impact of this minor temperature discrepancy 
is expected to be negligible compared to typical uncertainties in chemical concentrations. 

\begin{table}
    \centering
    \begin{tabular}{|c|}\hline
         \textbf{Kinetic parameters ($4^\circ C$)} \\\hline
         $k_{1a} = 6.2 \times 10^{-4}s^{-1}$  \\
         $k_{1b} = 5 \times 10^{-5} M$ \\
         $k_2 = 900 M^{-1}s^{-1}$ \\
         $k_{3a} = 100 M^{-2}s^{-2}$ \\
         $k_{3b} = 9.2 \times 10^{-5}s^{-1}$ \\
         $D_{MA} = 0.4 \times 10^{-5} cm^2s^{-1}$ \\
         $D_{I_2} = 0.6 \times 10^{-5}cm^2s^{-1}$ \\
         $D_{ClO_2} = 0.75 \times 10^{-5} cm^2s^{-1}$  \\
         $D_{I^-} = 0.7 \times 10^{-5}cm^2s^{-1}$ \\
         $D_{ClO_2^-} = 0.75 \times 10^{-5}cm^2s^{-1}$ \\\hline
         $\alpha = 10^{-15}$ \
         $k_4 = 10^8$ \
         $k_{-4} = 1$ \\\hline
    \end{tabular} \ \
%
    \begin{tabular}{|c|c|c|}\hline
     & Substance & Concentration  \\\hline
     Reservoir A & $\text{KIO}_3$ (Potassium iodate) & $1.88 \times 10^{-3} M$  \\
     & $\text{NaClO}_2$ (Sodium chlorite) & $2.0 \times 10^{-3} M$ \\
     & $\text{NaOH}$ (Sodium hydroxide) & $8.0 \times 10^{-3} M$ \\\hline
     Reservoir B & $\text{CH}_3\text{COOH}$ (Acetic acid) & $2.1 M$ \\\hline
     Hexagonal & KI (Potassium iodide) & $2.0 \times 10^{-3} M$ \\
     pattern & $\text{CH}_2(\text{COOH})_2$ (Malonic acid) & $2.25 \times 10^{-3} M$ \\\hline
     Striped & KI (Potassium iodide) & $2.0 \times 10^{-3} M$ \\
     pattern & $\text{CH}_2(\text{COOH})_2$ (Malonic acid) & $2.5 \times 10^{-3} M$ \\\hline
     Both & PVA (Polyvinylalcohol) & $1.5\, g / L$ \\\hline
\end{tabular}
    \caption{
    Left: Estimated kinetic parameters for the CIMA reaction at $4^\circ C$, as given in \cite{Rudovics1999}.\linebreak Right: Concentrations of substances fed into the reservoirs A and B of the two-sided open spatial reactor used in the chemical experiments  \cite{Rudovics1996} (shown in Figure \ref{fig:chemical_experiment}).}
    \label{tab:chemical_concentrations}
\end{table}

Next, the concentrations of $[\text{ClO}_2]$ and $[\text{KI}]$ in Reservoirs A and B need to be estimated based on the information provided in \cite{Rudovics1996} (see Table~\ref{tab:chemical_concentrations}, right). We assume that  $[\text{ClO}_2]$ is generated in Reservoir A from the reaction:
 \[
     \text{IO}_3^- + 6 \text{ClO}_2^- + 6 \text{H}^+   \rightarrow  6 \text{ClO}_2 + \text{I}^-  + 3 \text{H}_2\text{O},
 \]
 and thus set $[\text{ClO}_2] = [\text{NaClO}_2]_A = 2.0 \times 10^{-3} M$. To estimate the iodine concentration, we use the statement from \cite[p. 45]{Rudovics1996} that the iodine concentration follows that of iodide fed at the boundary of the gel. Thus, we set $[\text{I}_2] = 0.5 \times [\text{KI}]_B = 1.0 \times 10^{-3} M$. Additionally, we estimate $[\text{S}] = 1.5 \times 10^{-3}$ for 1.5\;$g/L$ of PVA, as done in \cite{Rudovics1999}. These computations are based on several assumptions and simplifications, and they are subject to measurement errors on several levels.   In particular, the concentrations inside the hydrogel block cannot be directly measured and can only be estimated from values on the boundary. Thus, all concentrations should be taken with a degree of caution.

 Using these estimates of the concentrations, rough estimates of model parameters for hexagonal and striped patterns can be derived from \eqref{eq:variable_change}. They are given in Table \ref{tab:parameter_estimates}. Note, however, that due to the uncertainties and inaccuracies discussed above these approximations may be off by as much as a factor of two or more, according to the opinion of experimentalists working in the field. In particular, we can see in Figure \ref{fig:results} that the Lengyel-Epstein model with these parameter values does not produce the observed patterns. 

The parameter estimates obtained by our algorithm in Table~\ref{tab:experiment_settings} are consistent with these crude experimentally-based predictions, when allowing for errors up to a factor 2-3. 
But what is more, our pattern-based numerical parameter estimation 
appears to have a significantly lower uncertainty and reliably allows to discern parameter regimes that lead to different pattern types present in the data. 



\section{Discussion}

In this paper, we presented an extension of the CIL statistical approach capable of robustly estimating the parameters of a pattern formation model from a single snapshot 
in the case of mixed mode patterns. We validated our method using real chemical data from the chlorite-iodite-malonic acid (CIMA) reaction. The pattern data comprised a low-resolution, noisy greyscale image of a small patch in the centre of a chemical reactor, without information regarding the actual concentration values of the observed chemical pattern or the initial conditions that lead to it. Notwithstanding this limitation, our method successfully fitted the Lengyel-Epstein reaction-diffusion model to the experimental pattern data, producing parameter estimates that were reasonably close to those predicted from chemical data. In the parameter identification process, we applied only basic pre-processing techniques, such as noise filtering, without leveraging any pattern-specific features of the data. Instead, we characterised the family of patterns using the distribution of distances.

The chemical experiment used to validate the proposed method closely replicates the challenges inherent in real biological systems. These include the inability to observe all model variables, the lack of information about the initial data that produced a given pattern, and the lack of absolute concentration values due to scaling. The proposed method’s ability to obtain robust parameter estimates despite these limitations makes it a valuable tool for various applications in biology. 

Moreover, our approach is not restricted to stationary patterns or specific pattern formation models. Although here we considered steady-state chemical patterns, the same methodology can be applied to spatio-temporal dynamics by including measurements taken at different time points during the transient phase, as was done in earlier works on the CIL approach for chaotic systems \cite{Haario2015,Springer2019,Springer2021}. Finally, the CIL approach is in principle not affected by the curse of dimensionality and scales well with respect to the number of parameters, in the same way as optimisation and MCMC sampling typically do, but we have not yet tested the approach on really high-dimensional problems. 



\subsection*{Acknowledgements}
We thank professor Philip Maini for suggesting the CIMA model as a real-world test problem for our approach.
This work was supported by the European Research Council (ERC) under the European Union’s Horizon 2020 research and innovation programme (synergy project PEPS, no. 101071786), and by the Deutsche Forschungsgemeinschaft (DFG) within the Collaborative Research Centre SFB1324 (B05) and within Germany’s Excellence Strategy EXC 2181/1 - 390900948 (the Heidelberg STRUCTURES Excellence Cluster). Additionally, this work was supported by the Research Council of Finland (RCoF) through the Flagship of advanced mathematics for sensing, imaging and modelling, decision number 358 944.

{\small
\bibliographystyle{acm}
\bibliography{cilexp_bibliography.bib}
}

\end{document}